\documentclass[10pt]{amsart}
\usepackage{amssymb,a4wide}
\usepackage{amsthm}
\usepackage{graphicx}
\usepackage[latin1]{inputenc}
\usepackage{tensor}
\usepackage{amscd}
\usepackage{url}
\usepackage{enumerate}

\newcommand{\weg}[1]{}

\theoremstyle{plain}
\newtheorem{thm}{Theorem}
\newtheorem*{thm*}{Theorem}

\theoremstyle{definition}

\theoremstyle{remark}
\newtheorem{rem}{Remark}

\begin{document}

\title[Quantum integrability for projectively equivalent metrics]{Quantum integrability for  the Beltrami-Laplace operators of projectively equivalent metrics of arbitrary signatures }
\author{Vladimir S. Matveev}
\address{Friedrich-Schiller-Universit\"at Jena, 
07737 Jena Germany  \ {\bf Email:}  vladimir.matveev@uni-jena.de}
\date{}

\begin{abstract} 
We generalize the result of \cite{ MT} to all signatures  
\end{abstract} 
 
\maketitle

 \vspace{-0.3cm}

\hfill {\small \it  Dedicated to Anatoly Timofeevich Fomenko  on his 75th birthday.}

 \vspace{0.3cm}

\section{Introduction}
Let $M$ be a smooth  manifold of dimension $n\ge 2$. We say that  two metrics $g$ and 
 $\bar g$  on this manifold are 
 \emph{projectively equivalent}, if each $g$-geodesic, after a proper reparameterization, is a $\bar g$-geodesic.  Theory of projectively equivalent metrics is a classical topic in differential geometry, already E. Beltrami  \cite{Beltrami} and T. Levi-Civita \cite{LC} did important contributions there.   In the  last two decades  a group of new methods    
 coming from integrable systems, see e.g.  \cite{MT1998, MT2000, Leeds2000, M2003,  cala},  and from Cartan geometry, see e.g. \cite{Eastwood,relativity, Iran}, appeared to be useful in this theory, 
  and made it possible to solve important open problems and named conjectures, see e.g. \cite{diffgeo,hyperbolic, Bryant, alone, mounoud}.

By \cite{MT1998,MT2003}  the existence of $\bar g$ projectively equivalent to $g$ allows one to construct a family $K^{(t)}_{ij}$ of Killing tensors of second degree for the metric $g$ (we will recall the formula and the definition later, in \S \ref{section:killing}, 
 following later publications, e.g. \cite{Bolsinov, hyperbolic, notes}. The family $K^{(t)}_{ij}$  is polynomial in $t$ of degree $n-1$ so it contains at most $n$ linearly independent Killing tensors).

In this paper we answer in Theorem \ref{thm:main} the following natural `quantization' question: {\it do the  corresponding second order differential operators commute?} 

There are of course many  possible  constructions of    differential operators of second order by   (0,2)-tensors, and, more generally, many different quantization approaches, see e.g. \cite[\S 6]{openproblems}.   We use the quantization procedure  of B. Carter   \cite[Equation (6.15)]{carter} and refer to  \cite{carter} and also to \cite{Chanu,Duval}  for an explanation why it  is natural in many aspects. The construction is as follows: to  a  tensor $K_{ij}$, we associate  an operator 
\begin{equation} \label{eq:quantisation}
\widehat K:C^\infty(M)\to C^\infty(M), \ \ \widehat K(f)= \nabla_i K^{ij} \nabla_j f.  
 \end{equation} 
Above and everywhere in the paper $\nabla $ is the Levi-Civita connection of $g$, we sum with respect to repeating indexes and raise the indexes of  $K$ by the metric $g$.

\begin{thm} \label{thm:main}  Assume  $g$ and $\bar g$ are projectively equivalent, let $K^{(t)}$ be the family of Killing tensors of second degree for  $g$  constructed with the help of $\bar g$.   Then, for any $t, s\in \mathbb{R} $,  the operators $\widehat K^{(t)}, \widehat K^{(s)}$   commute, that is  
$$
\widehat K^{(t)}  \widehat K^{(s)} - \widehat K^{(s)}  \widehat K^{(t)}=0.
$$
 \end{thm}   

Note that the Beltrami-Laplace operator $\Delta_g:= \nabla_i g^{ij} \nabla_j$ is a linear combination  of the operators of the family $\widehat K^{(t)}$, so all the operators $\widehat K^{(t)}$ commute also with $\Delta_g$. In fact, in the proof we go in the opposite direction:
we show first (combining \cite{carter,Duval} and \cite{KioMat2009})  that the operators  $\widehat K^{(t)}$  commute with  $\Delta_g$ and then use this to show that the operators  $\widehat K^{(t)}, \widehat K^{(s)}$ also commute mutually.

For Riemannian manifolds, Theorem \ref{thm:main} is known, it was announced in \cite{M} and the proof appeared in   \cite{MT}. The proof in the Riemannian case is based on direct  calculations in the coordinates in which the metrics admit the so-called Levi-Civita normal form. These coordinates exist (locally, in a neighborhood of almost every point),  if the (1,1)-tensor $G^i_j:= g^{is} \bar g_{sj}$ is semi-simple (at almost every point). This is always the case, for example, if one of the metrics is Riemannian. The proof from \cite{MT} can  be directly generalized to the pseudo-Riemannian metrics under the additional assumption that $G $ is semi-simple.

There are (many) examples of projectively equivalent metrics such that $G$ has  nontivial Jordan blocks; 
in this situation the proof and ideas of  \cite{MT} are  not  sufficient. Indeed, though also in  this  case  there exists a local description of  projectively equivalent metrics \cite{BM}, direct calculation  of the commutators of the operators $\widehat K^{(t)}$ and $\widehat K^{(s)}$
is  
 a complicated task     because of different combinatoric possibilities for the  number and the sizes of Jordan blocks and also because   the  description of \cite{BM} uses 
 a description of symmetric parallel (0,2)-tensors  from  \cite{boubel} which is quite nontrivial. For small dimensions it is possible though to prove Theorem \ref{thm:main} by direct calculations, in particual in dimension 2 it was done in \cite[\S 2.2.3]{pucacco}. 
 
 Our proof is based on another circle of ideas, it still uses the local description of \cite{BM} but replaces  local calculations by a trick which is based  on quite nontrivial results of different papers. We recall the necessary results in \S \ref{preliminary}.

  All objects in our paper are assumed to be sufficiently smooth. 
  
  We thank C. Chanu and V. Kiosak for useful discussions.

  \section{Basic facts about projectively equivalent metrics and Killing tensors  used in the proof} \label{preliminary}
  \subsection{Killing tensors  for projectively equivalent metrics and corresponding integrals.} \label{section:killing}
 Let $g$ and $\bar g$ be  two projectively equivalent metrics on the manifold $M$.   Let us 
  recall the construction of Killing tensors  $K^{(t)}_{ij}$ 
  of second degree for the metric $g$ by using  the  metric $\bar g$. 
  We consider the (1,1)-tensor $L$ given by the formula 
\begin{eqnarray}
L^i_j &:= & \left|\frac{\det(\bar g)}{\det(g)} \right|^\frac{1}{n+1}
 \bar g^{i\ell} g_{\ell j} \label{l}.
\end{eqnarray}
 Here $\bar g^{ij}$ is the contravariant metric
dual (= inverse, i.e., $\bar g^{is} \bar g_{sj}= \delta^i_j$)  to $\bar g$. 

Next, consider the family $S(t)$, $t\in \mathbb{R}$, of the 
$(1,1)$-tensors, where $\mbox{\rm Id}$ is the  (1,1)-tensor corresponding to the identity endomorphism, its components in the standard tensor notation are $\delta^i_j$. 
\begin{equation}\label{st}
 S(t):=   \operatorname{Comatrix}\left( t\ \mbox{\rm Id}-L\right).
 \end{equation}

Recall that the comatrix (or the  adjugate matrix) of a (1,1)-tensor is also a (1,1)-tensor.  Indeed, at  
 points where $t\not\in \operatorname{Spectrum}(L)$, 
 it is given by $$ \operatorname{Comatrix}\left(t\ \mbox{\rm Id}-L\right) =   \det\left(t\ \mbox{\rm Id}-L\right) \left(t\ \mbox{\rm Id}-L\right)^{-1}$$ and evidently  corresponds to a  (1,1)-tensor, and for each point the set of $t$ not lying in the spectrum of $L$ is everywhere  dense on the real line.   From the formula for the comatrix we see that the family \eqref{st} is polynomial in $t$  of degree $n-1$.

\begin{thm}[Essentially, \cite{MT1998}] \label{thm:killing} Let $g$ and $\bar g$ be projectively equivalent. Then, for 
 every  $t \in \mathbb{R}$ the tensor 
\begin{equation}\label{eq:killing}
K^{(t)}_{ij}:= g_{ir} S(t)^r_j
\end{equation}
is a Killing tensor for $g$. 
\end{thm}

In the coordinate-free notation the Killing tensor $K^{(t)}$ is given by $K^{(t)}( \xi , \nu )=  g( \xi, S(t)\nu)$. Since $L$ is $g$-selfadjoint, $S(t)$ is also self-adjoint so $K^{(t)}$ is  symmetric with respect to  the lower  indexes. Recall that a (symmetric with respect to  the lower  indexes) tensor $K_{ij}$ is Killing, if \begin{equation}\label{eq:killingequation} \nabla_{(i}K_{jk)}=0,\end{equation} 
where  the round  brackets denote the symmetrization. In our paper we do not use this equation, but use the geometric definition which we recall now:   a (0,2) symmetric tensor $K= K_{ij}$  is Killing, if and only if  the function 
$\tau\mapsto K(\gamma'(\tau), \gamma'(\tau)) $ is constant along every naturally  parameterized  $g$-geodesic $\gamma(\tau)$. In other words,  if the function  $   K(\gamma'(\tau), \gamma'(\tau)) $  is an integral of the geodesic flow of $g$. It is known, that the integrals corresponding to 
 the Killing tensors $K^{(t)}$ constructed above commute, let us  recall this statement:

\begin{thm} \label{thm:integrals} Let $g$ and $\bar g$ be projectively equivalent and $K^{(t)}$ be the Killing tensors for $g$ 
constructed  by \eqref{eq:killing}. Consider, for each $t\in \mathbb{R}$, 
 the function  $I_t:T^*M\to \mathbb{R}$  given by formula 
 \begin{equation} \label{eq:integral}
 I_t(x,p)= K^{(t)}_{rq}g^{iq} g^{ir} p_ip_j. 
 \end{equation}
 Here $(x,p)= (x^1,...,x^n,p_1,...,p_n)$ are    local  coordinates on $T^*M$: $x^i$ are local coordinates on $M$ and $p_i$ are, 
  for each $x$, the coordinates on $T^*_xM$ corresponding to the basis $\tfrac{\partial }{\partial x^i}$  on $T_xM$. 
  
  Then, for any $t, s\in \mathbb{R}$ the functions $I_t$, $I_s$ Poisson-commute with respect to the standard Poisson bracket on $T^*M$, that is:
  $$
  \sum_{i=1}^n \tfrac{\partial I_t}{\partial p_i} \tfrac{\partial I_s}{\partial x^i}-  \tfrac{\partial I_t}{\partial x^i} \tfrac{\partial I_s}{\partial p_i} =0.
  $$

\end{thm}

In the Riemannian signature, Theorem \ref{thm:integrals} is due to \cite{MT1998}. In all signatures, it was independently proved in \cite{Bolsinov,Topalov2001}.

  \subsection{ Difference between connections of projectively equivalent metrics} \label{section:d} 
  
  We consider the (1,1)-tensor $L$ constructed by projectively equivalent metrics $g$ and $\bar g$ by \eqref{l}. 
  As it was observed in \cite{Sinjukov}, see also \cite[Theorem 2]{Bolsinov},  it satisfies, for a certain 1-form $\lambda_i$,   the following equation: 
   \begin{equation} \label{basic} 
 \nabla_k L_{ij}= \lambda_i g_{jk} + \lambda_j  g_{ik}. 
 \end{equation} 
Here and later  we use $g$ for the covariant differentiations and for the tensor manipulations with indexes. By contracting \eqref{basic} with $g^{ij}$, we see that 
the 1-form $\lambda_i$  is the differential of the function $\lambda:= \tfrac{1}{2} \operatorname{trace}(L) = \tfrac{1}{2}  L^s_s.$

\begin{rem} The projectively-invariant form of this equation is due to \cite{Eastwood}, see also the survey  \cite{Iran} (and \cite{Bryant} for its two-dimensional version). It played essential role in many  recent  developments in the theory of projectively equivalent metrics  including the solutions of two problems explicitly stated by Sophus Lie \cite{Bryant,alone}, the  proof of the discrete version of the  projective 
Lichnerowciz conjecture \cite{fomenko70, zeghib} and the proof of the Lichnerowicz conjecture for metrics of Lorenzian signature \cite{BMR}. 
\end{rem} 

The 1-form $\lambda_i$ is closely related to the difference between the Levi-Civita 
connections of $\nabla= \left(\Gamma^i_{jk}\right) $ and $\bar \nabla= \left(\bar \Gamma^i_{jk}\right) $ (see e.g. \cite{Sinjukov} or \cite[\S 2.2]{KioMat2009}):  for the $1$-form 
\begin{equation}\label{phi}
\phi_i:= - L^s_i\lambda_s
\end{equation} 
we have \begin{equation} \label{cl} 
 \bar \Gamma_{jk}^i  - \Gamma_{jk}^i = \delta_k^i\phi_{j} + \delta_j^i\phi_{k}.    
   \end{equation}

From formulas (\ref{phi},\ref{cl}) we see that if $\lambda=\operatorname{trace}(L)$ has zero of order $k$ at a point $p\in M$, then at this point the connections coincide up to the order $k-1$. In particular, 
  for any tensor field $T$ the $(k-1)$st, and also lower order,  covariant derivatives of $T$ in $\nabla$ and $\bar \nabla$ coincide in $p$:
  $$
   \nabla_{i_1}\nabla_{i_2}...\nabla_{i_{k-1}} T \stackrel{\textrm{at $p$}}{=}  \bar \nabla_{i_1}\bar \nabla_{i_2}...\bar \nabla_{i_{k-1}}T.
  $$

Let us recall one more important property of projectively equivalent metrics: 

\begin{thm}[Folklore, e.g. Lemma 1 in  \cite{KioMat2009} or (12) in \cite{KioMat2010}]  \label{ricci} 
Let $g$ and $\bar g$ be projectively equivalent metrics and $L$ is as 
in \eqref{l}. Then, the Ricci curvature tensor $R_{ij}$ of $g$ commutes with $L$, in the sense
\begin{equation} \label{eq:comm}
R^i_s L^s_j- L^i_s R^s_j=0.
\end{equation}

\end{thm}

(For each $x\in M$ the 
formula \eqref{eq:comm}  is just the formula of the commutators of two endomorphisms of $T_xM$:
 the first is given by the Ricci tensor with one index raised, and the other it given by $L$).

\subsection{ Perturbing the metrics in the class of projectively equivalent metrics.} 
Let us now show that (for any $k$) one can perturb the metrics $g$ and $\bar g$ in the class of projectively equivalent metrics such that at a point they remain the same up to order $k$ and at another point the function $\lambda$ is constant up to order $k$.

 We say that two tensors or affine connections {\it coincide at a point $p$
 up to order $k$}, if  their difference is zero at $p$ and  in a local coordinate system 
   all  partial   derivatives   up to the  order $k$ of the components of their difference are zero at the point $p$. 
   This property does not depend on the choice of a coordinate system.
   
      In particular, a function is {\it  constant at $p$ up to order $k$ } 
       if all its partial derivatives up to  order  $k$  are zero at $p$.

\begin{thm}  \label{thm:order} Let $g$ and $\bar g$ be  projectively equivalent metrics and $L$ is as in \eqref{l}. Then, for each $k\in \mathbb{N}$ and for  almost any point $p\in M$  there exists  an arbitrary  small neighborhood $U$ containing $ p$,
a point $q\in U$ and  a pair of projectively  equivalent metrics $g' $ and $\bar g'$ on  $U$  (whose tensor \eqref{l} will be denoted by $L'$ and the function $ \tfrac{1}{2} \operatorname{trace}(L')$ will be denoted by $\lambda'$) 
 such that the following holds:  
\begin{itemize} \item[(A)] At  the  point $p$,   $g$ coincides with $g'$ and   $ \bar g$ coincides with $ \bar g'$    up to   order $ k$.
 \item[(B)] At the point $q$, $\lambda'$ is constant up to   order $k$.   
\end{itemize} 
\end{thm} 

 ``Almost every point''  means that the set of such  points  contains an open everywhere dense subset.  

   Theorem \ref{thm:order}  essentially follows from \cite{BM2011,BM}, let us explain this. We consider the points $p$ which are algebraically    generic  in the sense  of \cite[Def. 2.7]{Nijenhuis}: that is, there     exists a neighborhood $U\ni p$ such that at every point of the neighborhood        the number of different  eigenvalues of $L$ and   the number and the sizes of the Jordan blocks are the same 
   (of course  the eigenvalues are not necessary constant and usually depend on the point; by the implicit function theorem they are smooth functions near $p$).  

Take such a point. 
 Note that $\lambda$  is the half of the 
 sum of eigenvalues of $L$,
 counted with algebraic multiplicities. We need to find projectively equivalent 
 metrics $g'$ and $\bar g'$ 
such  
  that they  coincide  to order $k$ at  $p$ with $g$ and $\bar g$ 
   and such that  all eigenvalues of $L'$ are constant up to order $k$ in some point $q$. 
  
  By the Splitting-Gluing construction \cite[\S\S 1.1, 1.2]{BM2011},  it is sufficient to do this under the assumption that $L$ has one eigenvalue, or one pair of complex-conjugated eigenvalues.   
If the geometric multiplicity of an eigenvalue is greater than one, by  \cite[Proposition 1]{BM},  the eigenvalue is already a  constant, so we are done since $g'=g$ and $ \bar g= \bar g'$ are already as we want.

 Let us now consider the case when $L$ has one  real eigenvalue of   geometric multiplicity 1, or a pair of    nonreal  complex-conjugate eigenvalues of  geometric multiplicity 1. In this case,   
 the local structure of $g$ and $L$  near the point $p$ are described in some coordinate system. 
  There are 4 possible cases, the description was done in different papers, 
  let us give the precise  references where it  can be found.  
  
  If eigenvalue is real and its geometric  multiplicity is one (so  the ``splitted out''  manifold is one-dimensional), then 
the description is trivial  and was discussed e.g. 
 in  \cite[Example in \S 2.1]{BM2011} or \cite[Example 3 in \S 3.2.1]{relativity}.

   If $L$  has a pair of nonreal complex-conjugate eigenvalues of geometric multiplicity  $1$, then the description 
    was done in \cite[Theorem 2]{pucacco}, see also  \cite[Theorem A]{alone}.
    
     If $L$, at each point of $U$, 
   is conjugate to 
   a Jordan   block  with real eigenvalue, the description  is   in \cite[Theorem 4]{BM}.
   
    If $L$, at each point of $U$, 
   is conjugate to a pair of 
    Jordan   blocks   with complex-conjugated 
     eigenvalues, the description is done in   \cite[Theorem 5]{BM}.    
     
     In each of the above references, one sees that description is given by a formula and the only object we can choose is    the eigenvalue(s) of $L$: in the `real' case,    it  is  a function of one variable;  this function can be chosen arbitrary   
      (with exception that one may  not make it zero; though also this is allowed if we discuss not projectively equivalent metrics but `compatible'   in the terminology of \cite{BM}, pairs $(g, L)$).

       In the  `nonreal' case,    the eigenvalue   is a holomorphic function of one variable, again it can be chosen 
        arbitrary (again with exception that it is never zero) in the class of holomphic functions.  
  
In order to prove Theorem \ref{thm:order}, one  modifies  the  eigenvalue such that  at $p$ is coincides with the initial eigenvalue up to order $k$, and is constant up to order $k$ in some other point $q$. One can clearly do it for any function of one variable and for any holomorphic function of one complex variable.   

  \subsection{Carter's    condition.} We will need the following result:
  
 \begin{thm} \label{thm:carter}
Assume $K_{ij}$ is a Killing tensor for $g$ and $R_{ij}$ is the Ricci curvature tensor. Suppose, at the point $p\in M$, we have that up to order $k$
\begin{equation} \label{eq:carter} 
\nabla_i \left(R^i_s K^s_j- K^i_s R^s_j\right)= 0.
\end{equation} 
Then, the Beltrami-Lapalce operator $\Delta_g $ and the operator $\widehat K$ commute at the point $p$ up to order $k$, that is, for every function $f$ we have  
$$
\left(\Delta_g   \widehat K  -  \widehat K  \Delta_g  \right)f = 0 \ \ \textrm{at $p$  up to order $k$}
$$
\end{thm} 

Theorem above    is essentially due to B. Carter.  Indeed, from \cite[Equation (6.16)]{carter} it follows that 
 if $\nabla_i \left(R^i_s K^s_j- K^i_s R^s_j\right)$ is   zero \underline{at all points}, then $\Delta_g $ and   
 $\widehat K$ commute \underline{at all points}.  Careful analysis of the arguments 
 shows that   the proof of   Carter is valid also pointwise. Note that   only a sketch of the proof is given in \cite{carter}, and we recommend
  \cite[\S III(A)]{Duval} of   
 C. Duval and  G. Valent, from which  a more detailed proof can be extracted. More precisely, combining \cite[Equations (3.11) and (3.16)]{Duval} we obtain the above mentioned result of Carter.

\subsection{ If a Killing tensor vanishes up to a  sufficiently high  order  at one  point, then it is identically zero} 

\begin{thm} \label{jet}  Let $M$ be a connected manifold and $g$ be a metric of any signature on it. Assume $K$ is a Killing tensor of order $k$ (i.e., $K$ is a symmetric $(0,k)$ tensor satisfying the equation $\nabla_{(i} K_{i_1...i_k)}=0$). If $K$ vanishes up to order  $k$ at one point, then it vanishes identically on the whole manifold.  
\end{thm} 

This theorem follows from \cite{thompson} (see also \cite[\S 3]{kruglikov}).  We will need this theorem  for first and second degree  Killing tensors. Note that  for the first degree  Killing tensors (= Killing vectors, after raising the index),  Theorem  \ref{jet} can be obtained by the following geometric argument:    if a Killing vector field vanishes at a point $q$ up to order 1, then the flow of this vector field acts trivially on the tangent space to $q$. Since it commutes with the exponential mapping, the Killing vector field must be identically zero. For second  degree Killing tensors, the proof is based on the prolongation of the Killing equation which was essentially done in    \cite{wolf}. For all degree Killing tensors,  the prolongation  of Killing equation  was essentially done in \cite{thompson},
though formally this paper discusses special case of constant curvature metrics.    Indeed, for our goal  the  higher order terms of the prolongation are sufficient, and they do not depend on the curvature of the metric,   see e.g.  the discussion in \cite[\S 3]{kruglikov}). 

\section{ Proof of Theorem \ref{thm:main}.} 

We assume that $g$ and $\bar g$ are projectively equivalent metrics of any signature on $M^n$, $n\ge 2$. We
 consider $L$ given by \eqref{l},    the family 
$K^{(t)}$ of Killing tensors given  by  \eqref{eq:killing} and the corresponding differential  operators $\widehat K^{(t)}$. Combining Theorems \ref{ricci} and \ref{thm:carter}, we see that the operators commute with   $\Delta_g$.

Let us take any $t,s\in \mathbb{R}$ and consider the commutator 
$$
\widehat Q:= \widehat K^{(t)}  \widehat K^{(s)} - \widehat K^{(s)}  \widehat K^{(t)}.
$$
Our goal is to show that it vanishes; we will first show that  it is (linear) differential  operator of order  at most 2, i.e., that when we apply $\widehat Q$ to a function $f$ the  higher derivatives of $f$ vanish. This  step is well-known, see e.g.      \cite{carter} or \cite{Duval},  
let us shortly recall the arguments.

Clearly, $\widehat Q$ is a differential  operator of order at most 4, since both $\widehat K^{(t)}$ and 
$   \widehat K^{(s)} $ have order 2. One immediately sees though, 
that   the operators $\widehat K^{(t)}  \widehat K^{(s)}$  and $\widehat K^{(s)}  \widehat K^{(t)}$
have the same  symbols, so the 4th order terms  cancel when we subtract  one from the other. Thus,  the order of $\widehat Q$ is   at most 3. 
The third order terms vanish because the integrals corresponding to $K^{(s)}$ and $K^{(t)}$ commute by Theorem \ref{thm:integrals}. 
Indeed, direct calculations  show that the symbol of the commutator of two differential operators is the Poisson bracket of their symbols.

The proof   that the first and the  second order terms    vanish is  based on another (new) argument which will use all the results recalled  in \S \ref{preliminary}.

First observe  that  there exist a symmetric (2,0) tensor $Q^{ij}$ and the vector field $V^\ell$ such that 
 $$   
  \widehat Q =  \nabla_i Q^{ij} \nabla_j + V^\ell \nabla_\ell. 
  $$
 Indeed, the operator $
  \widehat Q$ does not have terms of zero order, since neither $ \widehat K^{(t)}$ nor $ \widehat K^{(s)}$ have such. One can collect all second order terms in $\nabla_i Q^{ij} \nabla_j$ and declare the rest as $V^\ell  \nabla_\ell$.

  Since $\Delta_g$ commutes with $\widehat K^{(t)} $ and $  \widehat K^{(s)}$,  it commutes with  $\widehat Q$. Then, $Q_{ij}$ is a Killing (0,2) tensor for $g$.

  It is sufficient to show, that $Q_{ij}$ vanishes at almost every point. It is sufficient to show this for almost every  $t$ and $s$.
    We take $s$ and $t$ such that the tensors $K^{(t)}, K^{(s)}$  are nondegenerate at some point. We will  work in a small neightborhood of this point, in each point of which the tensors $K^{(t)}, K^{(s)}$  are nondegenerate. Now we use Theorem  \ref{thm:order}: we first 
     take  a sufficiently big $k$ and then, for almost every point of $p$ of this neighborhood  
   consider the projectively equivalent      metrics $g'$ and $\bar g'$ satisfying conditions (A,B) from  Theorem \ref{thm:order}.  
    
     At the point $p$, the  metrics $g$ and $\bar g$ coincide with the metrics $g'$ and $\bar g'$,
      which implies that the  Killing  tensor   $Q'_{ij}$ (i.e., the analog  of the Killing tensor  $Q_{ij}$ constructed by $g'$ and $\bar g'$)     coincides with $Q_{ij}$ in $p$. 
     Let us show that, if $k$ is high enough,  at the point $q$ the Killing tensor $Q'_{ij}$ vanishes up to order 2. 
  
  At the point $q$, the 1-form $\lambda_i$ and therefore the 1-form $\phi_i$ (recalled in \S  \ref{section:d}) 
  vanishes up to (sufficiently high) order $k$. Then, at the point $q$, 
  the difference between Levi-Civita 
  connections $\nabla'$ of $g'$ and of $\bar g'$ vanishes up to order $k-1$, see \eqref{cl}. 
  Since the  Killing tensors   $K'^{(s)}$,  $K'^{(t)}$   are constructed by $g', \bar g'$ using algebraic formulas, 
   the covariant derivative in $\nabla'$ of    $K'^{(s)}$,  $K'^{(t)}$   vanishes at the point $q$ up to order $k-1$. 
  Then, up to the  order $k-1$, at the point $q$, 
  the Levi-Civita connection of the (contravariant)  metrics\footnote{As explicitly indicated, we view now  the Killing tensors as metrics: we first raise the indexes in   \eqref{eq:killing} by $g'$. The result is a nondegenerate symmetric $(2,0)$ tensor, we view it as a contravariant   metric. In order to obtain an usual metric, with lower   indexes, one needs to invert the matrix of  $\left(K'^{(t)}\right)^{ij}$.} $\left(K'^{(s)}\right)^{ij}$, $\left(K'^{(t)}\right)^{ij}$  coincide with $\nabla'$. 
  
  Then, at the point $q$, the Betrami-Laplace operators of the  the metrics  $K'^{(s)}$, $K'^{(t)}$ coincide with $\widehat K'^{(s)}$, $ \widehat K'^{(t)}$ up to order $k-2$. From the other side the Ricci tensor corresponding to the metric $K'^{(s)}$ commutes (in the sense of \eqref{eq:comm}) with  $K'^{(t)}$, up to the terms of order $k-3$, since it coincides up  to the terms of order $k-3$ with the Ricci tensor of $g'$ and it commutes with $L'$ and therefore with $S'(t)$. Then, the Carter condition  \eqref{eq:carter} is fulfilled up to order $k-4$.
   Then, the operators    $\hat K^{(t)}$ and $\hat K^{(s)}$
   commute at $q$ up to order $k-4$, which means that at $q$  we have  $Q'_{ij}=0$    up to  order $k-5$. If $k>7$, then this implies by Theorem \ref{thm:order}  that     $Q'_{ij}$ is identically zero, which means it vanishes at $p$, where it coincides with $Q_{ij}$. 
  Finally, $Q_{ij}= 0$ at $p$ and since $p$ was almost every point   $Q_{ij}\equiv 0$ on the whole manifold.

  \begin{rem} In fact the reader does  not need to follow the precise calculations of the necessary order above: it is clear that if $k$ is high enough then at the point $q$ the Levi-Civita connection of the 
   contravariant metric  corresponding to $K'^{(s)}$ (with upper indexes) coincides with that of $g$ up to a sufficiently high order and $K'^{(t)}$
    is parallel with respect to any of this connections up to a high order which means that the operators  
      $\hat K^{(t)}$ and $\hat K^{(s)}$
   commute at $p$ up to some high order and $Q'$ is zero up to a high order and is therefore identically zero.  
  \end{rem}

  But then $    
  \widehat Q =   V^\ell \nabla_\ell,  
  $
 since it commutes with $\Delta_g$,    $V^\ell$ is a Killing vector field. Using the same arguments, one shows that (for a perturbed metrics $g'$, $\bar g'$), $V'^\ell \equiv 0$, which implies that $V^\ell=0$ at $p$. Since this is fulfilled for almost all points $p$, we obtain   $V^{\ell}\equiv 0$.
Theorem \ref{thm:main} is proved.

  \section{Open problems}
  \subsection{Introducing potential} 
  
    We assume that $g$ and $\bar g$ are projectively equivalent metrics of any signature on $M^n$. We  consider the Killing tensors $K^{(t)}$ and the corresponding integrals $I_t$   from Theorem \ref{thm:integrals} and ask the following questions: \\
    
    \vspace{1ex}
     {\it  Can one add  functions $U^{(t)}:M\to \mathbb{R}$ to the integrals $I_t$  such that the results still Poisson-commute?  Do the corresponding differential operators, i.e.,  
     $\widehat K^{(t)} + U^{(t)}$,  still commute? } 
		
		\vspace{1ex}

		Of course it is interesting to get not one example of such functions (the trivial example $U^{(t)} =\operatorname{const}$ always exists)
		but construct all such examples, at least locally.

 If $L$ is semi-simple at almost every point (which is always the case if $g$ is Riemannian), the answer is positive,  which follows from the combination of results of \cite{k2,D}, see also \cite{ch}.   
  
  \subsection{Generalize the result for c-projectively equivalent metrics. }
   Theory of projectively equivalent metrics has  a natural analogue on Kähler manifolds: theory of c-projectively equivalent metrics.  Let us recall the basic definition:

Let $(M,g,J)$ be a  K\"ahler manifold of arbitrary signature of real dimension $2n\ge 4$. 
A regular curve $\gamma: \mathbb{R} \supseteq I \to M$ is called \emph{$J$-planar} if there exist functions $\alpha,\beta:I\rightarrow\mathbb{R}$ such that 
\begin{equation}
\label{eq:eq1} 
\nabla_{\dot \gamma(t)}\dot \gamma(t)=\alpha\dot \gamma(t)+\beta J(\dot \gamma(t))\mbox{ for all } t\in I,
\end{equation} 
where $\dot \gamma =\tfrac{d}{d t} \gamma$.

From the definition we see  immediately  that the property of $J$-planarity is independent of the parameterization of the curve, and that geodesics are $J$-planar curves.   
We also see that  $J$-planar curves form a much bigger family than the family of geodesics;  at every point and in every direction there exist 
infinitely many geometrically different  $J$-planar curves.

Two metrics $g$ and $\hat g$ of arbitrary signature  that are K\"ahler w.r.t the same complex structure $J$ are \emph{c-projectively equivalent} if any $J$-planar curve of $g$ is a $J$-planar curve of $\hat g$.  Actually, the condition that the metrics are K\"ahler with respect to the same complex structure is not essential; it is an easy exercise to show that if any $J$-planar curve of a K\"ahler structure $(g, J)$ is a $\hat J$-planar  curve of another K\"ahler structure 
$(\hat g, \hat J)$, then $\hat J=\pm J$.

C-projective equivalence was  introduced (under the name ``h-projective equivalence''or ``holomorphically projective correspondence'') 
 by T.  Otsuki and Y. Tashiro in  \cite{Otsuki,Tashiro}. Their motivation was to generalize the notion of 
 projective equivalence  to the K\"ahler situation.
Otsuki and Tashiro, see also  \cite[\S 6.2]{gover},  have shown that projective equivalence is not interesting in the  K\"ahler situation, since only simple examples are possible, and suggested c-projective equivalence as an  interesting object of study instead. This suggestion appeared to be very fruitful and between the 1960s and the 1970s, the theory of c-projectively equivalent metrics and c-projective transformations  was  one of the main research topics  in Japanese and Soviet (mostly Odessa and Kazan) differential geometry schools. Geometric structures  that are  equivalent to the existence of a   c-projective equivalent metric   were suggested independently in different branches of mathematics, see e.g. the introductions  of  \cite{R}  for a list and \cite{calderbank} for more detailed explanation on the relation to Hamiltonian 2-forms.

It appears that       many ideas and many results in the theory of projectively equivalent metrics
 have their counterparts in the c-projective setting.  For example, the use of 
 integrable systems  in the proof of the  Yano-Obata conjecture  \cite{Yano} about c-projective transformations 
  is very similar to that of in the Lichnerowicz conjecture \cite{diffgeo} for projective transformations. Compare also \cite{mettler, alone}.    See e.g. \cite[\S 1.2.]{BMR} for  one of the explanations. In particular, Theorems \ref{thm:killing} and  \ref{thm:integrals} have clear analogs:  by a c-projectively equivalent metric $\bar g$ 
 one can construct second degree Killing  tensors for $g$, and the  corresponding    integrals commute: see e.g. \cite[Proposition 5.14]{ECMN}, the result was initially obtained in \cite[Theorem 2]{Topalov2001}. 
 We ask the following question: can one generalize the result of the present  paper to c-projectively equivalent metrics? \\
 
 \vspace{1ex} 
 {\it Do the  differential operators corresponding to the Killing tensors from \cite[Proposition 5.14]{ECMN},  \cite[Theorem 2]{Topalov2001} commute?}
\vspace{1ex}

Also in the c-projective case, the Ricci tensor commutes with the analog of the tensor $L$. One can  do it  by the following 
tensor calculations which are similar to that of the proof of Theorem \ref{ricci}: take \cite[Equation (7)]{Ki} (which is the c-projective  analog of \cite[Equation (11)]{KioMat2010}), perturb the indexes by the  trivial permutation and by the 
   permutations    $ik\ell \mapsto  k\ell i$ and $ik\ell 
\mapsto  \ell ik$  and sum the results. 
We obtain  \cite[Equation (13)]{KioMat2010} (where  $a_{ij}$ corresponds to $L_{ij}$ in  our notation).    Contracting the obtained equation   with $g^{jk}$, we obtain an analog of   \eqref{eq:comm}, which implies by Theorem \ref{thm:carter} that 
the operators commute with the Beltrami-Laplace operator. 
Unfortunately,  the  rest of the proof  can not be directly generalized to  the c-projective case, since the analog of the  function 
$\lambda$ can not  be a constant up  to  high order  by  \cite[Corollary 3]{Ki}.   One can try to employ  \cite[Equation (3.11)]{Duval} for it, but we did not manage to overcome the technical   difficulties.

We do not have clear expectation how the answer would look: we tip that the operators do commute, but will not be suprised if    their commutators are first order differential operators corresponding to  Killing vector fields. We would like to recall  here that a  c-projectively equivalent metric allows one to construct Killing vector fields, see  e.g. \cite[\S  2]{BMR} and  \cite[\S 5.2]{ECMN}.

\end{document}